\documentclass{ifacconf}
\usepackage{xcolor}
\usepackage{enumerate}
\usepackage{amsmath, amssymb,mathrsfs}
\usepackage{natbib}            
\usepackage{graphicx}          
\makeatletter
\let\old@ssect\@ssect 
\makeatother

\usepackage{natbib}
\usepackage{hyperref}

\makeatletter
\def\@ssect#1#2#3#4#5#6{%
  \NR@gettitle{#6}
  \old@ssect{#1}{#2}{#3}{#4}{#5}{#6}
}
\makeatother
\usepackage{tikz,environ}
\usetikzlibrary{arrows,positioning,patterns,decorations.pathreplacing}

\tikzset{
    at xy split/.style 2 args={
        at={(#1,#2)}
    },
    a/.style={circle, draw=red},`'
    b/.style={rectangle, draw=blue}
}

\makeatletter
\newsavebox{\measure@tikzpicture}
\NewEnviron{scaletikzpicturetowidth}[1]{%
  \def\tikz@width{#1}%
  \begin{lrbox}{\measure@tikzpicture}%
  \BODY
  \end{lrbox}%
  \pgfmathparse{#1/\wd\measure@tikzpicture}%
  \BODY
}
\makeatother

\allowdisplaybreaks[3]

\begin{document}

\begin{frontmatter}

\title{A Minimax Optimal Controller for \\Positive Systems \thanksref{footnoteinfo}}

\thanks[footnoteinfo]{The authors are with the Department of Automatic Control and the ELLIIT Strategic Research Area at Lund University, Lund, Sweden. 
This work is partially funded by the Wallenberg AI, Autonomous Systems and Software Program (WASP) funded by the Knut and Alice Wallenberg Foundation, and the European Research Council (ERC) under the European Union's Horizon 2020 research and innovation programme under grant agreement No 834142  (ScalableControl).}

\author{Alba Gurpegui,} 
\author{Emma Tegling,} 
\author{Anders Rantzer}

\address{Automatic Control Department, Lund University, Sweden.  (e-mail: \{alba.gurpegui, emma.tegling, anders.rantzer\}@control.lth.se)}


\begin{abstract}                          
We present an explicit solution to the discrete-time Bellman equation for minimax optimal control of positive systems under unconstrained disturbances. The primary contribution of our result relies on deducing a bound for the disturbance penalty, which characterizes the existence of a finite solution to the problem class. Moreover, this constraint on the disturbance penalty reveals that, in scenarios where a solution is feasible, the problem converges to its equivalent minimization problem in the absence of disturbances.
\end{abstract}

\end{frontmatter}

\section{Introduction}

Minimax optimal control problems are widespread across control theory and engineering disciplines. They offer a robust methodology for formulating and addressing challenges characterized by competitive elements and uncertainties. Such problems are prevalent in areas including robust control, game theory, and multi-agent systems. Specifically, in our context, they are employed to design control systems that are robust to worst-case uncertainties and disturbances. Tackling these problems presents significant difficulties, particularly when we deal with large scale systems.

In this extended abstract we present a novel class of minimax optimal control problems, with positive dynamics, linear objective function, linear and homogeneous constraint on the control action, and unconstrained disturbances. We expanded upon the problem setting presented in~\cite{AlbaEmmaAnders}, removing the disturbance constraint.  Instead we introduce a bound for the disturbance penalty which characterizes the level of performance of the closed loop system and guarantees the existence of a finite solution. Furthermore, the condition on the disturbance penalty introduces a solution spectrum, highlighting the system's ability to maintain performance and stability in the face of disturbances. While dynamic programming is the basis we used to derive a solution in~\cite{AlbaEmmaAnders}, this extended work reveals that, if a solution is achievable, the problem aligns with its minimization class~\cite{mainArt}, whose explicit solution can be determined using linear programming. This solution uncovers that the optimal control policy, out of all potential policies, is linear in nature and characterized by a feedback matrix that mirrors the sparsity structure found in the problem's constraint matrix. This aspect is not just a theoretical curiosity, but has profound implications for the design of control systems. It enables the integration of structural controller constraints, particularly beneficial in large-scale systems where simplicity and scalability of implementation are paramount. In this context, the present work builds upon the approach of our prior study, optimizing sparse controllers without a priori constraints on linearity or sparsity and demonstrates that, for a broader class of problems where disturbances are unconstrained, among any potentially nonlinear and nonsparse controller, none can achieve a lower cost value than the one we present.

Another important characteristic of our problem class is the positive dynamics. Positive systems are characterized by the property that their state and output remain nonnegative for any nonnegative input and initial state. Classical references on this topic are~\cite{BermanBook} and~\cite{Luenberger}.
Extensive research has been conducted into expanding the scope of positive systems theory, exploring its natural extensions. These include positive systems with delays~\cite{Eibhara_steady_state_delay_interc_positive}, positive switched systems~\cite{Blanchini_switched_linear_positive}, and monotone systems~\cite{Smith1995MonotoneDS}. The application of positive systems theory spans a diverse array of fields. It has been instrumental in modeling dynamic systems in areas such as biology, ecology, physiology, and pharmacology~\cite{17tutorial, 19tutorial, 37tutorial, 38tutorial,39tutorial,44tutorial}, as well as in thermodynamics~\cite{Blanchini_switched_linear_positive, 38tutorial}. Its relevance extends to epidemiology~\cite{BLANCHINI_REV2, 40tutorial,53tutorial}, econometrics~\cite{55tutorial}, and even in more specific applications like filtering, charge routing networks, or power systems~\cite{6tutorial,7tutorial,10tutorial}.
One of the main advantages of positive systems is that stability can be verified using linear Lyapunov functions~\cite{Blanchini_lyapunov}, making this
class of systems more tractable in a large scale setting because of their computational scalability~\cite{EbiharaPeucelleArzelierTAC2017}. 

This extended abstract contains the explicit solution to our minimax optimal control problem setup and a sketch of the proof of its derivation. Here we indicate how this class of problems reduces to its minimization case, given the disturbance penalty condition we introduce in the main result of this manuscript. It is part of ongoing research to study the continuous setting and to compare it with past studies on $l_{1}$-controller synthesis for positive systems.
\subsection{Notation} Let $\mathbb{R}_{+}$ denote the set of nonnegative real numbers. The inequality $X > Y$ $(X \geq Y)$ means that all the elements of the matrix $(X-Y)$ are positive (nonnegative).  A matrix $X$ is called positive if all the elements of $X$ are nonnegative but at least one element is nonzero. The notation $\left | X \right |$ means element-wise absolute value, while $\mathbf{1}$ denotes a column vector with unit entries.
\subsection{Problem Setup}
    We present the optimal control problem of this extended abstract as a discrete-time, infinite-horizon, minimax optimal control problem with nonnegative cost and positive dynamics 
    \begin{align}\label{prob_setup}
    \underset{\mu}{\mathrm{inf}}\hspace{1mm}\underset{w}{\mathrm{\max}}&\sum_{t = 0}^{\infty} \left[s^{\top}x(t)+r^{\top}u(t)-\gamma^{\top}w(t) \right ]\\
    \mathrm{subject}& \hspace{1mm} \mathrm{to} \notag \\
    &x(t+1) = Ax(t)+Bu(t)+Fw(t), \notag \\
    &u(t)=\mu(x(t)) \hspace{1mm}; \hspace{2mm} x(0)=x_{0} \notag \\
    &\left | u \right |\leq Ex; \hspace{2mm} w \geq 0 \notag
    \end{align}
    where $x$ represents the $n$-dimensional vector of state variables, $u$ the $m$-dimensional control variable, $w$ the $l$-dimensional disturbance, $\mu$ is any, potentially nonlinear, control policy and $E$ prescribes the structure of the control action. The objective is to minimize the worst-case cost over all possible control strategies.
\section{Main Result}\label{mainresult}
We investigate the controller synthesis problem for~\eqref{prob_setup}. The main contribution is the derivation of the necessary condition~\eqref{gamma_ass}, for which a finite and explicit solution~\eqref{finite_value} is defined, the positive asymptotic stability is guaranteed and the level of performance $\gamma$ is satisfied.  
\begin{thm}
\label{MainThm} Let $A \in \mathbb{R}^{n\times n}$, $B = \left [ B_{1}~\cdots~ B_{m} \right ] \in \mathbb{R}^{n\times m }$, $F \in \mathbb{R}^{n \times l}_{+}$, $E = \left [ E_{1}^{\top} ~\cdots~E_{m}^{\top} \right ]^{\top} \in \mathbb{R}_{+}^{m \times n}$, $s \in \mathbb{R}^{n}$, $s>0$, $r \in \mathbb{R}^{m} $ and $\gamma \in \mathbb{R}^{l}_{+}$. Suppose that 
\begin{align}
    &A \geq \left | B \right |E \label{ass_A}\\
    &s > E ^{\top} \left | r \right |.\label{ass_s}
\end{align}
Then the optimal control problem~\eqref{prob_setup} has a finite value for every $x_{0} \in \mathbb{R}_{+}^{n}$ if and only if 
\begin{align}\label{gamma_ass}
    \gamma \geq F^{\top}p
\end{align}
where $p$ is obtained solving the linear program
\begin{align}\label{LProgram}
    &\mathrm{Maximize} \hspace{2mm} \mathbf{1}^{\top}p \hspace{1mm} \mathrm{over} \hspace{1mm} p \in \mathbb{R}^{n}_{+}, \hspace{1mm} \zeta \in \mathbb{R}^{m}_{+}\\
    &\mathrm{Subject} \hspace{1mm} \mathrm{to} \hspace{2mm} p\leq s+A^{\top}p-E^{\top}\zeta \notag \\
    &\hspace{15mm} -\zeta \leq r + B^{\top}p \leq \zeta . \notag
\end{align}
If this is true then~\eqref{prob_setup} has the minimum value $p^{\top}x_{0}$ with,
\begin{align}\label{finite_value}
    p & \leq s + A^{\top}p  -E^{\top} \left | r+B^{\top}p \right |^{\top}. 
\end{align}
Moreover, the control law $u(t) = -Kx(t)$, is optimal when 
\begin{align} \label{Kform}
    K := \begin{bmatrix}
\mathrm{sign}(r_{1}^{\top}+ p^{\top}B_{1})E_{1}\\ 
\vdots \\ 
\mathrm{sign}(r_{m}^{\top}+p^{\top}B_{m})E_{m}
\end{bmatrix}.
\end{align}
\end{thm}
\begin{rem}
\vspace{2mm}
The first condition~\eqref{ass_A} ensures the invariance of the positive orthant under the system dynamics.  
The second condition~\eqref{ass_s} will be needed when applying the assumptions of Lemma~\ref{LEMAA} in the Appendix to our objective function $g(x,u,w) = s^{T}x + r^{T}u -\gamma ^{T}w$.
\end{rem}
\begin{rem}
    The result in Theorem \ref{MainThm} is analogous for $w <0$ and $\gamma < F^{\top}p$ respectively.
\end{rem}
\begin{rem}
In \eqref{Kform}, it can be observed that the sparsity structure of the control gain $K$ is directly determined by the $E$ matrix. The sparsity of $E$ is in turn determined by the problem designer and may capture limitations in actuation and sensing. 
\end{rem}
\begin{pf} In this extended abstract we only provide an outline of the proof. The proof is based on Lemma~\ref{LEMAA} in the Appendix. First we reduce the general problem set up to our setting   
\begin{align*}
    f(x,u,w) &:= Ax+Bu+Fw \\
    g(x,u,w)&:= s^{\top}x+r^{\top}u - \gamma^{\top}w.
\end{align*}
where assumption~\eqref{ass_s} and~\eqref{ass_A} guarantee
\begin{align*}
\underset{w \geq 0}{\max}\left[g(x,u,w) \right] \geq 0.
\end{align*}
Next we use induction over $p^{\top}_{k}x=J_{k}(x)$ and  $p^{\top}x=J^{*}(x)$ for all $x$, to prove the equivalence between the recursive sequence $\left \{ p_{k}^{\top} \right \}_{k=0}^{\infty}$ and~\eqref{sequence_J} in Lemma~\ref{LEMAA}, and between the Bellman equation~\eqref{finite_value} and~\eqref{Bellman_EQ_g} in Lemma~\ref{LEMAA}. We use the equivalences in Lemma~\ref{LEMAA} to deduce the bound for the disturbance penalty~\eqref{gamma_ass} and prove implication $"\Rightarrow"$. Proving $"\Leftarrow"$ is direct, applying condition~\eqref{gamma_ass} to $J^{*}(x)$.   
\newline
\\
Finally, the expression for the optimal control policy is given by $u(t)=\mu(x(t))$, 
\begin{align*}
    \mu(x) &= \mathrm{arg}\min_{ | u|\leq Ex} 
    \left [ s^{\top}x+r^{\top}u-\gamma^{\top}w \right. \\
    &~~~~~~~~~~~~~~~~~~~~~~~~~~~~~~~~~~~~~~~\left.+ p^{\top}(Ax+Bu+Fw) \right ]\\
    &= \mathrm{arg}\min_{ | u|\leq Ex} \sum_{i=1}^{m}\left[ \left ( r_{i}^{\top}+p^{\top}B_{i} \right )u_{i} \right].
\end{align*}
and since for all $i = 1 \hspace{1mm}... \hspace{1mm} m$ the inequality $\left|u \right | \leq Ex$ restricts $u_{i}$ to the interval $\left [ -E_{i}x, E_{i}x \right ]$, the minimum is attained when $(r_{i}+p^{\top}B_{i})$ and $u_{i}$ have opposite signs. Thus, $u_{i} = -\mathrm{sign}(r_{i}+p^{\top}B_{i})E_{i}$ for all $i = 1 \hspace{1mm}... \hspace{1mm} m$.
\end{pf}
\section{Example. }\label{examples}
The optimal control problem~\eqref{prob_setup} admits sparsity constraints on the controller, making it particularly useful for large scale problems. Here, we consider a simple example of a double-tank process represented in Fig.~\ref{double_tank}, where the disturbance and the control action are equally characterized, because a large-scale example cannot be tractably represented in the text of this extended abstract. Instead we focus on the problem structure and the role of the constraints. The discretized double tank process dynamics in~\cite{double_tank} are given by: 
\begin{figure}[h]
   \centering
    \includegraphics[scale=0.12]{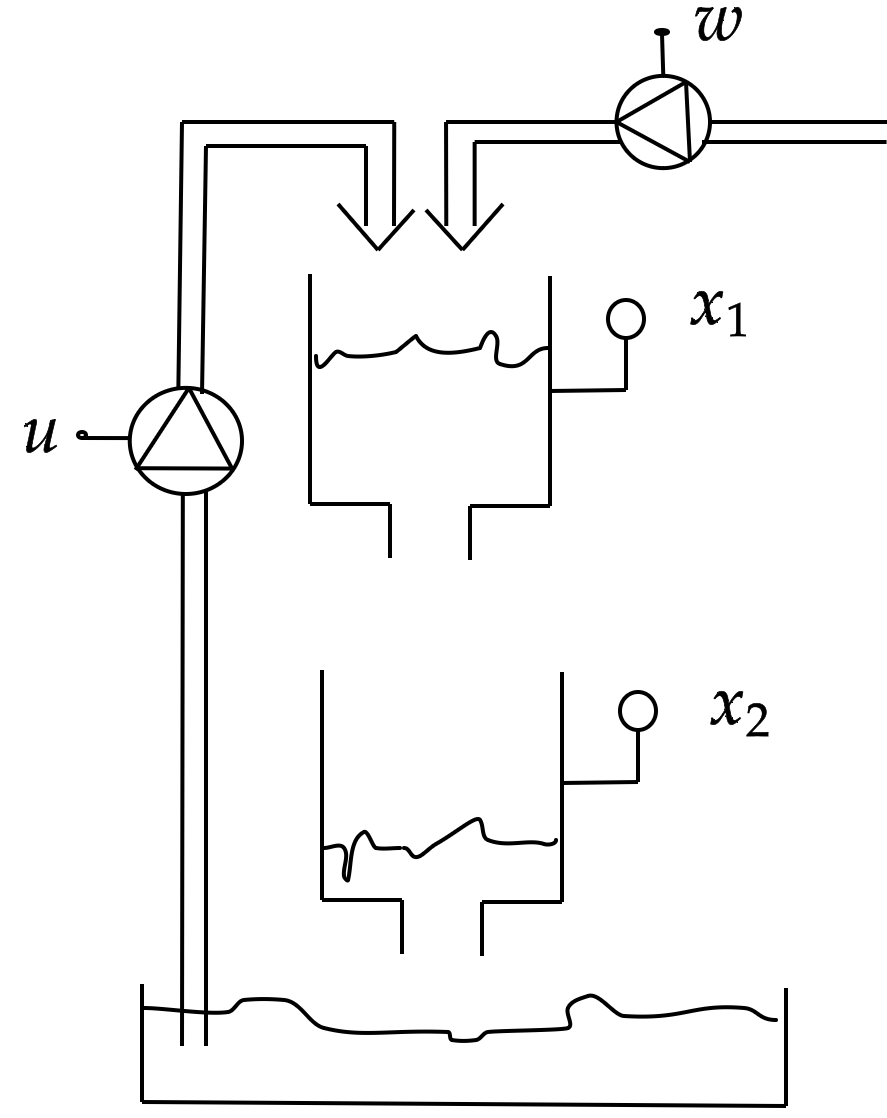}
    \caption{Double Tank process with disturbance $w$.}
    \label{double_tank}
\end{figure}
\begin{align*}
x(k+1) = A x(k) + B u(k) + F w(k)
\end{align*}
where 
\begin{align*}
    A = \begin{bmatrix}
0.9648 & 0\\ 
 0.0345 & 0.9648
\end{bmatrix}; \hspace{1mm} B = F=\begin{bmatrix}
0.0971 \\ 
 0.0017
\end{bmatrix}.
\end{align*}
Let $\gamma \geq F^{\top}p$, $E= \begin{bmatrix}
1 & 0
\end{bmatrix}$, $s= \begin{bmatrix}
1 & 1
\end{bmatrix}^{\top}$ and $r= 0.2 $. Now, we formulate the optimal control problem~\eqref{prob_setup} accordingly as:
\begin{align}\label{setup_doubletank}
    \underset{\mu}{\mathrm{inf}}\hspace{1mm}\underset{w}{\mathrm{\max}}&\sum_{t = 0}^{\infty} \left[s^{\top}x(t)+r^{\top}u(t)-\gamma^{\top}w(t) \right ]\\
    \mathrm{subject}& \hspace{1mm} \mathrm{to} \notag \\
    &x(t+1) = A x(t)+B u(t)+F w(t), \notag \\
    &u(t)=\mu(x(t)) \hspace{1mm}; \hspace{2mm} x(0)=x_{0} \notag \\
    &\left | u \right |\leq Ex; \hspace{2mm} w \geq 0 \notag
\end{align}
First we check that conditions~\eqref{ass_A} and~\eqref{ass_s} hold
\begin{align*}
    A =\begin{bmatrix}
0.9648 & 0\\ 
 0.0345 & 0.9648
\end{bmatrix} &\geq \begin{bmatrix}
0.0971 & 0\\ 
 0.0017 & 0
\end{bmatrix}=\left | B \right |E\\
s= \begin{bmatrix}
1 \\ 
1 
\end{bmatrix} &\geq \begin{bmatrix}
1 \\
0
\end{bmatrix}0.2 = E^{\top}r.
\end{align*}
Now, solving the linear program~\eqref{LProgram} we obtain 
\begin{align}\label{p_doubletank}
    p=\begin{bmatrix}
13.09 \\
28.41
\end{bmatrix}; \zeta =1.52.
\end{align}
Finally, from Theorem~\ref{MainThm} we know that the problem~\eqref{setup_doubletank} has a solution if and only if $\gamma \geq F^{\top}p \approx 1.32$. For this value of $\gamma$ the optimal solution is~\eqref{p_doubletank} and its respective feedback controller matrix is $K=E$.

\bibliography{ifacconf}
\section*{Appendix}
The proof of the following results are presented in \cite{AlbaEmmaAnders}.
We define, a general discrete time, infinite horizon, minimax optimal control problem with discrete cost function and constraints as
\begin{align}
    \label{gen_setup}
    \underset{\mu}{\mathrm{inf}}\hspace{1mm}\underset{w}{\mathrm{\max}}&\sum_{t = 0}^{\infty} g(x(t), u(t), w(t)) \\
    \mathrm{subject}& \hspace{1mm} \mathrm{to} \notag \\
    &x(t+1) = f(x(t), u(t), w(t)), \notag \\
    & x(t)\in X; \hspace{2mm} x(0)=x_{0}\hspace{1mm}; \hspace{2mm} u(t)=\mu(x(t)) \notag \\
    &u(t)\in U(x(t)); \hspace{2mm} w(t)\in W(x(t)) \notag
\end{align}
where $f: X \times U \times W \rightarrow X$, $x$ represents the vector of $n$-dimensional state variables, $u$ the $m$-dimensional control variable and $w$ the $q$-dimensional disturbance. 

\begin{lem} \label{LEMAA}
Suppose 
\begin{align*}
    \underset{w\in W(x)}{\mathrm{max}}\left[g(x,u,w) \right] \geq 0
\end{align*}
$\forall x \in X$, $\forall u \in U(x)$. Then, the following statements are equivalent. 
\begin{enumerate}[(i)]
    \item  The general optimal control problem in~\eqref{gen_setup} has a finite value for every $x_{0} \in \mathbb{R}_{+}^{n}$. 
    \vspace{1mm}
    
    \item The recursive sequence $\left \{ J_{k} \right \}_{k=0}^{\infty}$ with $J_{0} =0$ and
    \begin{align}\label{sequence_J}
         J_{k}(x)=\underset{u }{\mathrm{\min} }\hspace{1mm}\underset{w }{\mathrm{\max}}\left [ g(x,u,w)+J_{k-1}(f(x,u,w)) \right ]
    \end{align}
    has a finite limit $\forall x \in X$.
    \vspace{1mm}
    
    \item The Bellman equation  
    \begin{align}\label{Bellman_EQ_g}
        J^{*}(x)&=\underset{u}{\mathrm{min}} \hspace{1mm} \underset{w }{\mathrm{max}}\left [ g(x,u,w) + J^{*}(f(x,u,w))\right ] 
    \end{align} 
    has nonnegative solution $J^{*}(x)$, $\forall x \in X$.
\end{enumerate}
\end{lem}
\end{document}